\documentclass[12pt]{article}%
\usepackage{amsfonts}
\usepackage{sw20bams}
\usepackage{amsmath}
\usepackage{amssymb}
\usepackage{graphicx}%
\providecommand{\U}[1]{\protect\rule{.1in}{.1in}}
\begin{document}

\title{Moments of Maximum: Segment of AR(1)}
\author{Steven Finch}
\date{August 12, 2019}
\maketitle

\begin{abstract}
Let $X_{t}$ denote a stationary first-order autoregressive process. \ Consider
five contiguous observations (in time $t$) of the series (e.g., $X_{1}$,
\ldots, $X_{5}$). \ Let $M$ denote the maximum of these. \ Let $\rho$ be the
lag-one serial correlation, which satisfies $\left\vert \rho\right\vert <1$.
\ For what value of $\rho$ is $\mathbb{E}(M)$ maximized? \ How does
$\mathbb{V}(M)$ behave for increasing $\rho$? \ Answers to these questions lie
in Afonja (1972), suitably decoded.

\end{abstract}

\footnotetext{Copyright \copyright \ 2019 by Steven R. Finch. All rights
reserved.}Assume that $(X_{1},\ldots,X_{\ell})$ is $\ell$-multivariate
normally distributed with vector mean zero and covariance matrix%
\[
R=\left(
\begin{array}
[c]{ccccc}%
1 & \rho_{12} & \rho_{13} & \ldots & \rho_{1\ell}\\
\rho_{21} & 1 & \rho_{23} & \ldots & \rho_{2\ell}\\
\rho_{31} & \rho_{32} & 1 & \ldots & \rho_{3\ell}\\
\vdots & \vdots & \vdots & \ddots & \vdots\\
\rho_{\ell1} & \rho_{\ell2} & \rho_{\ell3} & \ldots & 1
\end{array}
\right)  .
\]
In particular, all variances are one and the correlation between $X_{i}$ and
$X_{j}$ is $\rho_{ij}=\rho_{ji}$. \ Define%
\[
\Phi_{\ell}(R)=\mathbb{P}\left\{  X_{i}\geq0\text{ for all }i=1,\ldots
,\ell\right\}
\]
and $M=\max\{X_{1},\ldots,X_{\ell}\}$. \ Afonja \cite{Afj-heu} proved that%
\[
\mathbb{E}\left(  M\right)  =%
%TCIMACRO{\dsum \limits_{i=1}^{\ell}}%
%BeginExpansion
{\displaystyle\sum\limits_{i=1}^{\ell}}
%EndExpansion%
%TCIMACRO{\dsum \limits_{j\neq i}}%
%BeginExpansion
{\displaystyle\sum\limits_{j\neq i}}
%EndExpansion
\sqrt{\frac{1-\rho_{ij}}{4\pi}}\cdot\Phi_{\ell-2}(R_{i,j}),
\]%
\[
\mathbb{E}\left(  M^{2}\right)  =1+%
%TCIMACRO{\dsum \limits_{i=1}^{\ell}}%
%BeginExpansion
{\displaystyle\sum\limits_{i=1}^{\ell}}
%EndExpansion%
%TCIMACRO{\dsum \limits_{j\neq i}}%
%BeginExpansion
{\displaystyle\sum\limits_{j\neq i}}
%EndExpansion%
%TCIMACRO{\dsum \limits_{\substack{k\neq i\\k\neq j}}}%
%BeginExpansion
{\displaystyle\sum\limits_{\substack{k\neq i\\k\neq j}}}
%EndExpansion
\frac{1-\rho_{ij}}{4\pi}\cdot\frac{1+\rho_{ij}-\rho_{ik}-\rho_{jk}}%
{\sqrt{4(1-\rho_{ij})(1-\rho_{ik})-(1-\rho_{ij}-\rho_{ik}+\rho_{jk})^{2}}%
}\cdot\Phi_{\ell-3}(R_{i,,jk})
\]
where the matrices $R_{i,j}$ and $R_{i,jk}$ require specification. \ As a
preliminary step, let%
\begin{align*}
r_{i,jk}  &  =\left\{
\begin{array}
[c]{lll}%
\text{correlation between }X_{i}-X_{j}\text{ and }X_{i}-X_{k} &  & \text{if
}j\neq i\text{ and }k\neq i,\\
\text{correlation between }X_{i}-X_{j}\text{ and }X_{i} &  & \text{if }j\neq
i\text{ and }k=i,\\
\text{correlation between }X_{i}\text{ and }X_{i}-X_{k} &  & \text{if
}j=i\text{ and }k\neq i,\\
1 &  & \text{if }j=i\text{ and }k=i
\end{array}
\right. \\
&  =\left\{
\begin{array}
[c]{lll}%
\dfrac{1-\rho_{ij}-\rho_{ik}+\rho_{jk}}{\sqrt{2(1-\rho_{ij})}\sqrt
{2(1-\rho_{ik})}} &  & \text{if }j\neq i\text{ and }k\neq i,\\
\dfrac{1-\rho_{ij}}{\sqrt{2(1-\rho_{ij})}} &  & \text{if }j\neq i\text{ and
}k=i,\\
\dfrac{1-\rho_{ik}}{\sqrt{2(1-\rho_{ik})}} &  & \text{if }j=i\text{ and }k\neq
i\\
1 &  & \text{if }j=i\text{ and }k=i
\end{array}
\right. \\
&  =\left\{
\begin{array}
[c]{lll}%
\dfrac{1-\rho_{ij}-\rho_{ik}+\rho_{jk}}{\sqrt{4(1-\rho_{ij})(1-\rho_{ik})}} &
& \text{if }j\neq i\text{ and }k\neq i,\\
\sqrt{\dfrac{1-\rho_{ij}}{2}} &  & \text{if }j\neq i\text{ and }k=i,\\
\sqrt{\dfrac{1-\rho_{ik}}{2}} &  & \text{if }j=i\text{ and }k\neq i,\\
1 &  & \text{if }j=i\text{ and }k=i
\end{array}
\right.
\end{align*}
and
\[
R_{i}=\left(
\begin{array}
[c]{ccc}%
r_{i,11} & \ldots & r_{i,1\ell}\\
\vdots &  & \vdots\\
r_{i,\ell1} & \ldots & r_{i,\ell\ell}%
\end{array}
\right)  .
\]
Clearly $r_{i,jk}=r_{i,kj}$ always. \ The matrix $R_{i}$ will not be needed in
ensuing sections -- we have included it for completeness' sake only -- it
would, however, come into play if means were nonconstant.

\section{Partial Correlations}

Set $\ell=4$ for simplicity. \ Fix $1\leq i<j\leq4$. \ There exists a unique
ordered pair $m<n$ such that $(i,j,m,n\}=\{1,2,3,4\}$. \ Define a $3\times3$
matrix%
\[
P=\left(
\begin{array}
[c]{ccc}%
1 & r_{i,mn} & r_{i,jm}\\
r_{i,mn} & 1 & r_{i,jn}\\
r_{i,jm} & r_{i,jn} & 1
\end{array}
\right)
\]
which captures all correlations among $X_{i}-X_{m}$, $X_{i}-X_{n}$,
$X_{i}-X_{j}$. \ Let $P_{ab}$ be the cofactor of the element $p_{ab}$ in the
expansion of the determinant of $P$. \ The partial correlation $r_{i,mn.j}$
between $X_{i}-X_{m}$ and $X_{i}-X_{n}$, given $X_{i}-X_{j}$, is prescribed by%
\[
r_{i,mn.j}=-\frac{P_{12}}{\sqrt{P_{11}P_{22}}}=\frac{r_{i,mn}-r_{i,jm}%
r_{i,jn}}{\sqrt{\left(  1-r_{i,jm}^{2}\right)  \left(  1-r_{i,jn}^{2}\right)
}}.
\]
The notation here \cite{Afj-heu} differs from that used
elsewhere\ \cite{Prk-heu, BS-heu, Fi1-heu}. \ In words, $r_{i,mn.j}$ measures
the linear dependence of $X_{i}-X_{m}$ and $X_{i}-X_{n}$ in which the
influence of $X_{i}-X_{j}$ is removed. \ Define finally a $2\times2$ matrix%
\[
R_{i,j}=\left(
\begin{array}
[c]{cc}%
1 & r_{i,mn.j}\\
r_{i,mn.j} & 1
\end{array}
\right)
\]
as was to be done.

Set now $\ell=5$. \ Fix $1\leq i<j\leq5$. \ There exists a unique ordered
triple $m<n<o$ such that $(i,j,m,n,o\}=\{1,2,3,4,5\}$. \ The preceding
discussion extends naturally, supplementing the case $(m,n)$ by additional
possible cases $(m,o)$ and $(n,o)$. \ Define finally a $3\times3$ matrix%
\[
R_{i,j}=\left(
\begin{array}
[c]{ccc}%
1 & r_{i,mn.j} & r_{i,mo.j}\\
r_{i,mn.j} & 1 & r_{i,no.j}\\
r_{i,mo.j} & r_{i,no.j} & 1
\end{array}
\right)  .
\]
We could go on for larger $\ell$, but this is all that is needed for our purposes.

Again, set $\ell=5$. \ Fix $1\leq i<j<k\leq5$. \ There exists a unique ordered
pair $m<n$ such that $(i,j,k,m,n\}=\{1,2,3,4,5\}$. \ Define a $4\times4$
matrix%
\[
Q=\left(
\begin{array}
[c]{cccc}%
1 & r_{i,mn} & r_{i,jm} & r_{i,km}\\
r_{i,mn} & 1 & r_{i,jn} & r_{i,kn}\\
r_{i,jm} & r_{i,jn} & 1 & r_{i,jk}\\
r_{i,km} & r_{i,kn} & r_{i,jk} & 1
\end{array}
\right)
\]
which captures all correlations among $X_{i}-X_{m}$, $X_{i}-X_{n}$,
$X_{i}-X_{j}$, $X_{i}-X_{k}$. \ Let $Q_{ab}$ be the cofactor of the element
$q_{ab}$ in the expansion of the determinant of $Q$. \ The partial correlation
$r_{i,mn.jk}$ between $X_{i}-X_{m}$ and $X_{i}-X_{n}$, given $X_{i}-X_{j}$ and
$X_{i}-X_{k}$, is prescribed by%
\begin{align*}
r_{i,mn.jk}  &  =-\frac{Q_{12}}{\sqrt{Q_{11}Q_{22}}}\\
&  =\frac{r_{i,mn}-r_{i,jm}r_{i,jn}-r_{i,km}r_{i,kn}+r_{i,km}r_{i,jn}%
r_{i,jk}+r_{i,jm}r_{i,kn}r_{i,jk}-r_{i,mn}r_{i,jk}^{2}}{\sqrt{\left(
1-r_{i,jm}^{2}-r_{i,km}^{2}+2r_{i,jm}r_{i,km}r_{i,jk}-r_{i,jk}^{2}\right)
\left(  1-r_{i,jn}^{2}-r_{i,kn}^{2}+2r_{i,jn}r_{i,kn}r_{i,jk}-r_{i,jk}%
^{2}\right)  }}.
\end{align*}
In words, $r_{i,mn.jk}$ measures the linear dependence of $X_{i}-X_{m}$ and
$X_{i}-X_{n}$ in which the influence of $X_{i}-X_{j}$ and $X_{i}-X_{k}$ is
removed. \ Define finally a $2\times2$ matrix%
\[
R_{i,jk}=\left(
\begin{array}
[c]{cc}%
1 & r_{i,mn.jk}\\
r_{i,mn.jk} & 1
\end{array}
\right)
\]
as was to be done.

Set now $\ell=6$. \ Fix $1\leq i<j<k\leq6$. \ There exists a unique ordered
triple $m<n<o$ such that $(i,j,k,m,n,o\}=\{1,2,3,4,5,6\}$. \ The preceding
discussion extends naturally, supplementing the case $(m,n)$ by additional
possible cases $(m,o)$ and $(n,o)$. \ Define finally a $3\times3$ matrix%
\[
R_{i,jk}=\left(
\begin{array}
[c]{ccc}%
1 & r_{i,mn.jk} & r_{i,mo.jk}\\
r_{i,mn.jk} & 1 & r_{i,no.jk}\\
r_{i,mo.jk} & r_{i,no.jk} & 1
\end{array}
\right)  .
\]
We could go on for larger $\ell$, but this is all that is needed for our purposes.

\section{Small Segments}

For convenience, define%
\begin{align*}
h(x,y,z)  &  =\frac{1-x}{4\pi}\cdot\frac{1+x-y-z}{\sqrt
{4(1-x)(1-y)-(1-x-y+z)^{2}}}\\
&  =\frac{1-x}{4\pi}\cdot\frac{1+x-y-z}{\sqrt{4(1-x)(1-z)-(1-x-z+y)^{2}}}\\
&  =\frac{1-x}{4\pi}\cdot\frac{1+x-y-z}{\sqrt
{(1-x)(3+x)-y(2+y)-z(2+z)+2(xy+xz+yz)}}.
\end{align*}
The latter expression, while more cumbersome, exhibits symmetry in $y$, $z$. \ 

If $\ell=2$, then $\Phi_{\ell-2}(R_{i,j})=1$ and $\Phi_{\ell-3}(R_{i,,jk})=0$.
\ We have%
\[%
\begin{array}
[c]{ccc}%
\mathbb{E}\left(  M\right)  =2\sqrt{\dfrac{1-\rho_{12}}{4\pi}}=\sqrt
{\dfrac{1-\rho_{12}}{\pi}}, &  & \mathbb{E}\left(  M^{2}\right)  =1.
\end{array}
\]

If $\ell=3$, then $\Phi_{\ell-2}(R_{i,j})=\frac{1}{2}$ and $\Phi_{\ell
-3}(R_{i,,jk})=1$. \ We have%
\[
\mathbb{E}\left(  M\right)  =\frac{1}{\sqrt{4\pi}}\left(  \sqrt{1-\rho_{12}%
}+\sqrt{1-\rho_{13}}+\sqrt{1-\rho_{23}}\right)  ,
\]%
\[
\mathbb{E}\left(  M^{2}\right)  =1+2h(\rho_{12},\rho_{13},\rho_{23}%
)+2h(\rho_{13},\rho_{12},\rho_{23})+2h(\rho_{23},\rho_{12},\rho_{13}).
\]
In formula (3.6)\ for $\mathbb{E}\left(  M\right)  $ in \cite{Afj-heu},
$\sqrt{(\pi)}$ should be replaced by $\sqrt{(2\pi)}$. \ 

If $\ell=4$, then
\[%
\begin{array}
[c]{ccc}%
\Phi_{\ell-2}(R_{i,j})=\frac{1}{4}+\frac{1}{2\pi}\arcsin(r_{i,mn.j}), &  &
\Phi_{\ell-3}(R_{i,,jk})=\frac{1}{2}.
\end{array}
\]
In general, $r_{i,mn.j}\neq r_{j,mn.i}$ and thus symmetry fails for
$\mathbb{E}\left(  M\right)  $. \ We have%
\begin{align*}
\mathbb{E}\left(  M\right)   &  =\tfrac{1}{\sqrt{4\pi}}\left[  \sqrt
{1-\rho_{12}}\cdot\left(  \tfrac{1}{4}+\tfrac{1}{2\pi}\arcsin(r_{1,34.2}%
)\right)  +\sqrt{1-\rho_{12}}\cdot\left(  \tfrac{1}{4}+\tfrac{1}{2\pi}%
\arcsin(r_{2,34.1})\right)  \right. \\
&  +\sqrt{1-\rho_{13}}\cdot\left(  \tfrac{1}{4}+\tfrac{1}{2\pi}\arcsin
(r_{1,24.3})\right)  +\sqrt{1-\rho_{13}}\cdot\left(  \tfrac{1}{4}+\tfrac
{1}{2\pi}\arcsin(r_{3,24.1})\right) \\
&  +\sqrt{1-\rho_{14}}\cdot\left(  \tfrac{1}{4}+\tfrac{1}{2\pi}\arcsin
(r_{1,23.4})\right)  +\sqrt{1-\rho_{14}}\cdot\left(  \tfrac{1}{4}+\tfrac
{1}{2\pi}\arcsin(r_{4,23.1})\right) \\
&  +\sqrt{1-\rho_{23}}\cdot\left(  \tfrac{1}{4}+\tfrac{1}{2\pi}\arcsin
(r_{2,14.3})\right)  +\sqrt{1-\rho_{23}}\cdot\left(  \tfrac{1}{4}+\tfrac
{1}{2\pi}\arcsin(r_{3,14.2})\right) \\
&  +\sqrt{1-\rho_{24}}\cdot\left(  \tfrac{1}{4}+\tfrac{1}{2\pi}\arcsin
(r_{2,13.4})\right)  +\sqrt{1-\rho_{24}}\cdot\left(  \tfrac{1}{4}+\tfrac
{1}{2\pi}\arcsin(r_{4,13.2})\right) \\
&  \left.  +\sqrt{1-\rho_{34}}\cdot\left(  \tfrac{1}{4}+\tfrac{1}{2\pi}%
\arcsin(r_{3,12.4})\right)  +\sqrt{1-\rho_{34}}\cdot\left(  \tfrac{1}%
{4}+\tfrac{1}{2\pi}\arcsin(r_{4,12.3})\right)  \right]  ,
\end{align*}%
\begin{align*}
\mathbb{E}\left(  M^{2}\right)   &  =1+h(\rho_{12},\rho_{13},\rho_{23}%
)+h(\rho_{12},\rho_{14},\rho_{24})+h(\rho_{13},\rho_{12},\rho_{23}%
)+h(\rho_{13},\rho_{14},\rho_{34})\\
&  +h(\rho_{14},\rho_{12},\rho_{24})+h(\rho_{14},\rho_{13},\rho_{34}%
)+h(\rho_{23},\rho_{12},\rho_{13})+h(\rho_{23},\rho_{24},\rho_{34})\\
&  +h(\rho_{24},\rho_{12},\rho_{14})+h(\rho_{24},\rho_{23},\rho_{34}%
)+h(\rho_{34},\rho_{13},\rho_{14})+h(\rho_{34},\rho_{23},\rho_{24}).
\end{align*}
In formula (3.7)\ for $\mathbb{E}\left(  M\right)  $ in \cite{Afj-heu}, a
factor $1/\sqrt{2\pi}$ should be inserted in front of the summation. \ 

If $\ell=5$, then%
\[
\Phi_{\ell-2}(R_{i,j})=\tfrac{1}{2}-\tfrac{1}{4\pi}\arccos(r_{i,mn.j}%
)-\tfrac{1}{4\pi}\arccos(r_{i,mo.j})-\tfrac{1}{4\pi}\arccos(r_{i,no.j}),
\]%
\[
\Phi_{\ell-3}(R_{i,,jk})=\tfrac{1}{4}+\tfrac{1}{2\pi}\arcsin(r_{i,mn.jk}).
\]
Symmetry now fails for both $\mathbb{E}\left(  M\right)  $ and $\mathbb{E}%
\left(  M^{2}\right)  $. \ We have%
\begin{align*}
\mathbb{E}\left(  M\right)   &  =\tfrac{1}{\sqrt{4\pi}}\left[  \sqrt
{1-\rho_{12}}\cdot\left(  \tfrac{1}{2}-\tfrac{1}{4\pi}\arccos(r_{1,34.2}%
)-\tfrac{1}{4\pi}\arccos(r_{1,35.2})-\tfrac{1}{4\pi}\arccos(r_{1,45.2}%
)\right)  \right. \\
&  +\sqrt{1-\rho_{12}}\cdot\left(  \tfrac{1}{2}-\tfrac{1}{4\pi}\arccos
(r_{2,34.1})-\tfrac{1}{4\pi}\arccos(r_{2,35.1})-\tfrac{1}{4\pi}\arccos
(r_{2,45.1})\right) \\
&  +\sqrt{1-\rho_{13}}\cdot\left(  \tfrac{1}{2}-\tfrac{1}{4\pi}\arccos
(r_{1,24.3})-\tfrac{1}{4\pi}\arccos(r_{1,25.3})-\tfrac{1}{4\pi}\arccos
(r_{1,45.3})\right) \\
&  \left.  +\sqrt{1-\rho_{13}}\cdot\left(  \tfrac{1}{2}-\tfrac{1}{4\pi}%
\arccos(r_{3,24.1})-\tfrac{1}{4\pi}\arccos(r_{3,25.1})-\tfrac{1}{4\pi}%
\arccos(r_{3,45.1})\right)  +\cdots\right]  ,
\end{align*}%
\begin{align*}
\mathbb{E}\left(  M^{2}\right)   &  =1+h(\rho_{12},\rho_{13},\rho_{23}%
)\cdot\left(  \tfrac{1}{4}+\tfrac{1}{2\pi}\arcsin(r_{1,45.23})\right)
+h(\rho_{12},\rho_{23},\rho_{13})\cdot\left(  \tfrac{1}{4}+\tfrac{1}{2\pi
}\arcsin(r_{2,45.13})\right) \\
&  +h(\rho_{12},\rho_{14},\rho_{24})\cdot\left(  \tfrac{1}{4}+\tfrac{1}{2\pi
}\arcsin(r_{1,35.24})\right)  +h(\rho_{12},\rho_{24},\rho_{14})\cdot\left(
\tfrac{1}{4}+\tfrac{1}{2\pi}\arcsin(r_{2,35.14})\right) \\
&  +h(\rho_{12},\rho_{15},\rho_{25})\cdot\left(  \tfrac{1}{4}+\tfrac{1}{2\pi
}\arcsin(r_{1,34.25})\right)  +h(\rho_{12},\rho_{25},\rho_{15})\cdot\left(
\tfrac{1}{4}+\tfrac{1}{2\pi}\arcsin(r_{2,34.15})\right)  +\cdots,
\end{align*}
a total of $20$ terms and $61$ terms, respectively.

If $\ell=6$, then $\mathbb{E}\left(  M\right)  $ contains non-elementary
functions which require numerical integration (beyond our present scope). \ In
contrast,
\begin{align*}
\mathbb{E}\left(  M^{2}\right)   &  =1+h(\rho_{12},\rho_{13},\rho_{23}%
)\cdot\left(  \tfrac{1}{2}-\tfrac{1}{4\pi}\arccos(r_{1,45.23})-\tfrac{1}{4\pi
}\arccos(r_{1,46.23})-\tfrac{1}{4\pi}\arccos(r_{1,56.23})\right) \\
&  +h(\rho_{12},\rho_{23},\rho_{13})\cdot\left(  \tfrac{1}{2}-\tfrac{1}{4\pi
}\arccos(r_{2,45.13})-\tfrac{1}{4\pi}\arccos(r_{2,46.13})-\tfrac{1}{4\pi
}\arccos(r_{2,56.13})\right) \\
&  +h(\rho_{12},\rho_{14},\rho_{24})\cdot\left(  \tfrac{1}{2}-\tfrac{1}{4\pi
}\arccos(r_{1,35.24})-\tfrac{1}{4\pi}\arccos(r_{1,36.24})-\tfrac{1}{4\pi
}\arccos(r_{1,56.24})\right) \\
&  +h(\rho_{12},\rho_{24},\rho_{14})\cdot\left(  \tfrac{1}{2}-\tfrac{1}{4\pi
}\arccos(r_{2,35.14})-\tfrac{1}{4\pi}\arccos(r_{2,36.14})-\tfrac{1}{4\pi
}\arccos(r_{2,56.14})\right) \\
&  +h(\rho_{12},\rho_{15},\rho_{25})\cdot\left(  \tfrac{1}{2}-\tfrac{1}{4\pi
}\arccos(r_{1,34.25})-\tfrac{1}{4\pi}\arccos(r_{1,36.25})-\tfrac{1}{4\pi
}\arccos(r_{1,46.25})\right) \\
&  +h(\rho_{12},\rho_{25},\rho_{15})\cdot\left(  \tfrac{1}{2}-\tfrac{1}{4\pi
}\arccos(r_{2,34.15})-\tfrac{1}{4\pi}\arccos(r_{2,36.15})-\tfrac{1}{4\pi
}\arccos(r_{2,46.15})\right) \\
&  +h(\rho_{12},\rho_{16},\rho_{26})\cdot\left(  \tfrac{1}{2}-\tfrac{1}{4\pi
}\arccos(r_{1,34.26})-\tfrac{1}{4\pi}\arccos(r_{1,35.26})-\tfrac{1}{4\pi
}\arccos(r_{1,45.26})\right) \\
&  +h(\rho_{12},\rho_{26},\rho_{16})\cdot\left(  \tfrac{1}{2}-\tfrac{1}{4\pi
}\arccos(r_{2,34.16})-\tfrac{1}{4\pi}\arccos(r_{2,35.16})-\tfrac{1}{4\pi
}\arccos(r_{2,45.16})\right)  +\cdots,
\end{align*}
a total of $121$ terms. \ In formula (3.9)\ for $\mathbb{E}\left(
M^{2}\right)  $ in \cite{Afj-heu}, a constant term $1$ should be inserted in
front of the first summation; further, the last summation should be taken over
both $k\neq i$ and $k\neq j$ (not merely $k\neq i$).

\section{Time Series}

Consider a discrete-time stationary first-order autoregressive process%
\[%
\begin{array}
[c]{ccccc}%
X_{t}=\rho\,X_{t-1}+\sqrt{1-\rho^{2}}\cdot\varepsilon_{t}, &  & -\infty
<t<\infty, &  & |\rho|<1
\end{array}
\]
where $\varepsilon_{t}$ is $N(0,1)$ white noise. The $\ell\times\ell$
covariance matrix $R$ has $ij^{\text{th}}$ element
\[
\rho_{ij}=\rho^{\left\vert j-i\right\vert }%
\]
which leads to certain simplifications. \ Let us make the reliance of $M$ on
$\ell$ explicit. \ We have%
\[
\mathbb{E}\left(  M_{2}\right)  =\sqrt{\dfrac{1-\rho}{\pi}},
\]%
\[
\mathbb{E}\left(  M_{3}\right)  =\sqrt{\dfrac{1-\rho}{\pi}}+\sqrt
{\dfrac{1-\rho^{2}}{4\pi}},
\]%
\begin{align*}
&  \mathbb{E}\left(  M_{4}\right)  =\tfrac{1}{4\pi}\sqrt{\tfrac{1-\rho}{\pi}%
}\left[  \pi+2\arcsin\left(  1-\tfrac{2}{3-\rho}\right)  \right]  +\tfrac
{1}{2\pi}\sqrt{\tfrac{1-\rho}{\pi}}\left[  \pi+2\arcsin\left(  \tfrac
{1+2\rho-\rho^{2}}{\sqrt{(3-\rho)\left(  3+\rho+\rho^{2}-\rho^{3}\right)  }%
}\right)  \right]  \\
&  +\tfrac{1}{2\pi}\sqrt{\tfrac{1-\rho^{2}}{\pi}}\left[  \pi+2\arcsin\left(
\tfrac{(1-\rho)^{2}}{\sqrt{(3-\rho)\left(  3+\rho+\rho^{2}-\rho^{3}\right)  }%
}\right)  \right]  +\tfrac{1}{4\pi}\sqrt{\tfrac{1-\rho^{3}}{\pi}}\left[
\pi+2\arcsin\left(  1-\tfrac{2}{3+\rho+\rho^{2}-\rho^{3}}\right)  \right]
\end{align*}
but $\mathbb{E}\left(  M_{5}\right)  $ is too lengthy to record here. \ In the
limit as $\rho\rightarrow0$, we obtain%
\[
\frac{1}{\sqrt{\pi}},\;\;\frac{3}{2\sqrt{\pi}},\;\;\frac{3}{\sqrt{\pi}}\left[
1-\frac{1}{\pi}\operatorname*{arcsec}(3)\right]  ,\;\;\frac{5}{\sqrt{\pi}%
}\left[  1-\frac{3}{2\pi}\operatorname*{arcsec}(3)\right]
\]
for $\ell=2,3,4,5$ and these are consistent with well-known values
\cite{Fi0-heu} corresponding to independent $X_{t}$. \ Figure 1 displays
$\mathbb{E}\left(  M_{\ell}\right)  $ as functions of $\rho$. \ The left-hand
endpoint is at $(-1,\sqrt{2/\pi})$ and $\sqrt{2/\pi}$ is unsurprisingly the
mean of a standard half-normal distribution. \ The right-hand endpoint is at
$(1,0)$. \ Associated with $\ell=3,4,5$ are maximum points with $\rho$ equal
to%
\[
\frac{1-\sqrt{5}}{2}=-0.6180339887498948482045868...,
\]%
\[%
\begin{array}
[c]{ccc}%
-0.4973597615161907364022217..., &  & -0.4336476843162656141275672...
\end{array}
\]
respectively. \ Closed-form expressions for the latter two quantities remain open.

We also have%
\[
\mathbb{E}\left(  M_{3}^{2}\right)  =1+\frac{(1-\rho)\sqrt{(3-\rho)(1+\rho)}%
}{2\pi},
\]%
\[
\mathbb{E}\left(  M_{4}^{2}\right)  =1+\tfrac{3+\sqrt{(3-\rho)\left(
3+\rho+\rho^{2}-\rho^{3}\right)  }+\rho\left[  1-2\rho-2\rho^{2}-\rho^{3}%
+\rho^{4}-\rho\sqrt{(3-\rho)\left(  3+\rho+\rho^{2}-\rho^{3}\right)  }\right]
}{2\pi\sqrt{(1+\rho)\left(  3+\rho+\rho^{2}-\rho^{3}\right)  }}%
\]
but $\mathbb{E}\left(  M_{5}^{2}\right)  $ and $\mathbb{E}\left(  M_{6}%
^{2}\right)  $ are too lengthy to record here. \ In the limit as
$\rho\rightarrow0$, we obtain%
\[
1+\frac{\sqrt{3}}{2\pi},\;\;1+\frac{\sqrt{3}}{\pi},\;\;1+\frac{5\sqrt{3}}%
{2\pi}\left[  1-\frac{1}{\pi}\operatorname*{arcsec}(4)\right]  ,\;\;1+\frac
{5\sqrt{3}}{\pi}\left[  1-\frac{3}{2\pi}\operatorname*{arcsec}(4)\right]
\]
for $\ell=3,4,5,6$ and these again are consistent with well-known values
\cite{Fi0-heu}. \ Associated with $\ell=3,4,5,6$ are maximum points with
$\rho$ equal to
\[
1-\sqrt{2}=-0.4142135623730950488016887,
\]%
\[%
\begin{array}
[c]{ccc}%
-0.3879232988398265768779440..., &  & -0.3599267104829689555367968...
\end{array}
\]%
\[
-0.3406053067160525788737944...
\]
respectively. \ Closed-form expressions for the latter three quantities remain
open. \ Figure 2 displays $\mathbb{V}\left(  M_{\ell}\right)  =\mathbb{E}%
\left(  M_{\ell}^{2}\right)  -\mathbb{E}\left(  M_{\ell}\right)  ^{2}$ as
functions of $\rho$. \ The left-hand endpoint is at $(-1,1-2/\pi)$; the
right-hand endpoint is at $(1,1)$. \ Unlike $\mathbb{E}\left(  M_{\ell
}\right)  $ or $\mathbb{E}\left(  M_{\ell}^{2}\right)  $, the variance is
strictly increasing throughout the interval. \ An intuitive reason for such
behavior would be good to establish someday.

\section{Proof of Revision}

Our general formula for $\mathbb{E}\left(  M^{2}\right)  $ looks somewhat
different from that presented by Afonja \cite{Afj-heu}. \ To demonstrate the
equivalence of the two formulas, it suffices to prove that if $j\neq i$,
$k\neq i$ and $k\neq j$, then%
\[
\frac{1}{2\pi}r_{i,ji}\cdot\frac{r_{i,ki}-r_{i,jk}r_{i,ji}}{\sqrt
{1-r_{i,jk}^{2}}}=h\left(  \rho_{ij},\rho_{ik},\rho_{jk}\right)  .
\]
The left-hand side is equal to%
\begin{align*}
&  \dfrac{1}{2\pi}\sqrt{\dfrac{1-\rho_{ij}}{2}}\cdot\frac{\sqrt{\dfrac
{1-\rho_{ik}}{2}}-\dfrac{1-\rho_{ij}-\rho_{ik}+\rho_{jk}}{\sqrt{4(1-\rho
_{ij})(1-\rho_{ik})}}\sqrt{\dfrac{1-\rho_{ij}}{2}}}{\sqrt{1-\dfrac
{(1-\rho_{ij}-\rho_{ik}+\rho_{jk})^{2}}{4(1-\rho_{ij})(1-\rho_{ik})}}}\\
&  =\dfrac{1}{4\pi}\sqrt{1-\rho_{ij}}\cdot\frac{\sqrt{1-\rho_{ik}}%
-\dfrac{1-\rho_{ij}-\rho_{ik}+\rho_{jk}}{\sqrt{4(1-\rho_{ik})}}}%
{\sqrt{1-\dfrac{(1-\rho_{ij}-\rho_{ik}+\rho_{jk})^{2}}{4(1-\rho_{ij}%
)(1-\rho_{ik})}}}\\
&  =\dfrac{1}{4\pi}\sqrt{1-\rho_{ij}}\cdot\frac{2(1-\rho_{ik})-(1-\rho
_{ij}-\rho_{ik}+\rho_{jk})}{\sqrt{4(1-\rho_{ik})-\dfrac{(1-\rho_{ij}-\rho
_{ik}+\rho_{jk})^{2}}{1-\rho_{ij}}}}\\
&  =\dfrac{1}{4\pi}(1-\rho_{ij})\cdot\frac{1+\rho_{ij}-\rho_{ik}-\rho_{jk}%
}{\sqrt{4(1-\rho_{ij})(1-\rho_{ik})-(1-\rho_{ij}-\rho_{ik}+\rho_{jk})^{2}}}%
\end{align*}
which is the right-hand side, as was to be shown.

\section{Proof from First Principles}

An exercise in \cite{DN-heu} suggests that formulas for $\mathbb{E}\left(
M_{2}\right)  $ and $\mathbb{V}\left(  M_{2}\right)  $ should be derived from%
\[
\max\left\{  X_{1},X_{2}\right\}  =\frac{1}{2}(X_{1}+X_{2})+\frac{1}%
{2}\left\vert X_{1}-X_{2}\right\vert .
\]
It is instructive to similarly prove our formula for $\mathbb{E}\left(
M_{3}\right)  $, using instead%
\begin{align*}
\max\left\{  X_{1},X_{2},X_{2}\right\}   &  =\max\left\{  \max\left\{
X_{1},X_{2}\right\}  ,\max\left\{  X_{2},X_{3}\right\}  \right\} \\
&  =\frac{1}{2}\left(  \max\left\{  X_{1},X_{2}\right\}  +\max\left\{
X_{2},X_{3}\right\}  \right)  +\frac{1}{2}\left\vert \max\left\{  X_{1}%
,X_{2}\right\}  -\max\left\{  X_{2},X_{3}\right\}  \right\vert \\
&  =\frac{1}{4}\left(  (X_{1}+X_{2})+(X_{2}+X_{3})+\left\vert X_{1}%
-X_{2}\right\vert +\left\vert X_{2}-X_{3}\right\vert \right) \\
&  +\frac{1}{4}\left\vert (X_{1}+X_{2})-(X_{2}+X_{3})+\left\vert X_{1}%
-X_{2}\right\vert -\left\vert X_{2}-X_{3}\right\vert \right\vert .
\end{align*}
Define $Y=X_{1}-X_{2}$ and $Z=X_{3}-X_{2}$. \ Clearly $(Y,Z)$ is bivariate
normally distributed with vector mean zero and covariance matrix%
\[
\left(
\begin{array}
[c]{cc}%
2-2\rho_{12} & \rho_{13}-\rho_{12}-\rho_{23}+1\\
\rho_{13}-\rho_{12}-\rho_{23}+1 & 2-2\rho_{23}%
\end{array}
\right)  =\left(
\begin{array}
[c]{cc}%
\sigma_{y}^{2} & \xi\,\sigma_{y}\sigma_{z}\\
\xi\,\sigma_{y}\sigma_{z} & \sigma_{z}^{2}%
\end{array}
\right)  .
\]
Also%
\[
\frac{1}{4}\mathbb{\,E\,}\left\vert Y\right\vert =\frac{\sigma_{y}}{4}%
\sqrt{\frac{2}{\pi}}=\sqrt{\dfrac{1-\rho_{12}}{4\pi}},
\]%
\[
\frac{1}{4}\mathbb{\,E\,}\left\vert Z\right\vert =\frac{\sigma_{z}}{4}%
\sqrt{\frac{2}{\pi}}=\sqrt{\dfrac{1-\rho_{23}}{4\pi}}.
\]
The four integrals (depending on signs of $Y$ and $Z$) underlying%
\[
\frac{1}{4}\mathbb{\,E\,}\left\vert (Y+\left\vert Y\right\vert )-(Z+\left\vert
Z\right\vert )\right\vert =\sqrt{\dfrac{1-\rho_{13}}{4\pi}}%
\]
can all be evaluated (however tediously). \ Because%
\[
X_{1}-X_{3}=(X_{1}-X_{2})-(X_{3}-X_{2})=Y-Z
\]
we suspect that a more elegant proof ought to be available. \ Ideas on
bridging this gap would be welcome.

In more detail, letting%
\[
f(y,z)=\frac{1}{2\pi\sqrt{1-\xi^{2}}\,\sigma_{y}\sigma_{z}}\exp\left[
-\frac{1}{2\left(  1-\xi^{2}\right)  }\left(  \frac{y^{2}}{\sigma_{y}^{2}%
}-\frac{2\xi yz}{\sigma_{y}\sigma_{z}}+\frac{z^{2}}{\sigma_{z}^{2}}\right)
\right]
\]
denote the bivariate normal density, we obtain%
\[
\int\limits_{0}^{\infty}\int\limits_{0}^{\infty}2\left\vert y-z\right\vert
f(y,z)dy\,dz=\frac{1}{\sqrt{2\pi}}\left[  -(1-\xi)(\sigma_{y}+\sigma
_{z})+2\sqrt{\sigma_{y}^{2}-2\xi\sigma_{y}\sigma_{z}+\sigma_{z}^{2}}\right]
\]
when $Y>0$ and $Z>0$;
\[
\int\limits_{0}^{\infty}\int\limits_{-\infty}^{0}2\,z\,f(y,z)dy\,dz=\frac
{(1-\xi)\sigma_{z}}{\sqrt{2\pi}}%
\]
when $Y<0$ and $Z>0$;%
\[
\int\limits_{-\infty}^{0}\int\limits_{0}^{\infty}2\,y\,f(y,z)dy\,dz=\frac
{(1-\xi)\sigma_{y}}{\sqrt{2\pi}}%
\]
when $Y>0$ and $Z<0$; and $0$ when $Y<0$ and $Z<0$. \ Adding these
contributions and dividing by $4$, we verify%
\begin{align*}
\frac{1}{2\sqrt{2\pi}}\sqrt{\sigma_{y}^{2}-2\xi\sigma_{y}\sigma_{z}+\sigma
_{z}^{2}} &  =\frac{1}{2\sqrt{\pi}}\sqrt{(1-\rho_{12})-(\rho_{13}-\rho
_{12}-\rho_{23}+1)+(1-\rho_{23})}\\
&  =\sqrt{\dfrac{1-\rho_{13}}{4\pi}}%
\end{align*}
as was desired.

Calculating the variance of $M_{3}$ from first principles\ has not been
attempted. \ The variance of the median ($50\%$-tile) is also of interest,
appearing explicitly in \cite{Fi5-heu} for $\ell=3$ but under the assumption
of independence. \ 

An alternative probability density-based derivation of $\mathbb{E}\left(
M_{2}\right)  $ and $\mathbb{V}\left(  M_{2}\right)  $ can be found in
\cite{Htr-heu, NK-heu}. \ See also \cite{Fi2-heu} for the expected range of a
normal sample, \cite{Fi3-heu} for the expected \textit{absolute} maximum, and
\cite{Fi4-heu} for other aspects of AR(1).

\section{Large Segments}

Assuming $\rho_{ij}$ depends only on $\left\vert j-i\right\vert =d$, Berman
\cite{Be-heu, Pi-heu, Ja-heu} proved that if either%
\[%
\begin{array}
[c]{ccccc}%
\lim\limits_{d\rightarrow\infty}\rho(d)\cdot\ln(d)=0 &  & \text{or} &  &
%TCIMACRO{\dsum \limits_{d=1}^{\infty}}%
%BeginExpansion
{\displaystyle\sum\limits_{d=1}^{\infty}}
%EndExpansion
\rho(d)^{2}<\infty,
\end{array}
\]
then%
\[
\lim\limits_{\ell\rightarrow\infty}P\left\{  \sqrt{2\ln(\ell)}(M_{\ell
}-a_{\ell})\leq x\right\}  =\exp\left(  -e^{-x}\right)
\]
where%
\[
a_{\ell}=\sqrt{2\ln(\ell)}-\frac{1}{2}\frac{\ln(\ln(\ell))+\ln(4\pi)}%
{\sqrt{2\ln(\ell)}}.
\]
Further, the two hypotheses on $\rho(d)$ cannot be significantly weakened.
\ This theorem clearly applies for a first-order autoregressive process,
although we note that $a_{\ell}$ does not incorporate lag-one correlation
$\rho$ at all. \ A more precise asymptotic result might do so. \ 

Other relevant works in the literature include \cite{R1-heu, R2-heu, WM-heu,
AG-heu, W1-heu, W2-heu, W3-heu, CBS-heu, Na-heu}. \ In particular, Figure 2 of
\cite{AG-heu} depicts the density of AR(1) maximum for $\ell=5$ and
$\rho=-9/10$, $-8/10$, \ldots, $8/10$, $9/10$. \ 

\section{Acknowledgements}

Raymond Kan \cite{KR-heu} symbolically evaluated the integrals in Section 5,
at my request, for the special case $\sigma_{y}=\sigma_{z}=1$ using
Mathematica. \ Enrique del Castillo \cite{CBS-heu} identified several
typographical errors in \cite{WM-heu} and provided R\ code for numerically
approximating the first two moments\ of $M_{\ell}$. \ I am grateful to both
individuals for their kindness!

%

%TCIMACRO{\FRAME{ftbpFU}{6.3005in}{4.1594in}{0pt}{\Qcb{ Interior maximum points
%exist for $\ell\geq3$, but not for $\ell=2.$}}{}{afonjam0.eps}%
%width 6.3005in;  height 4.1594in;  depth 0pt;  original-width 10.6692in;
%original-height 7.0197in;  cropleft "0";  croptop "1";  cropright "1";
%cropbottom "0";  filename 'afonjam0.eps';file-properties "XNPEU";}}}%
%BeginExpansion
\begin{figure}[ptb]%
\centering
\includegraphics[
height=4.1594in,
width=6.3005in
]%
{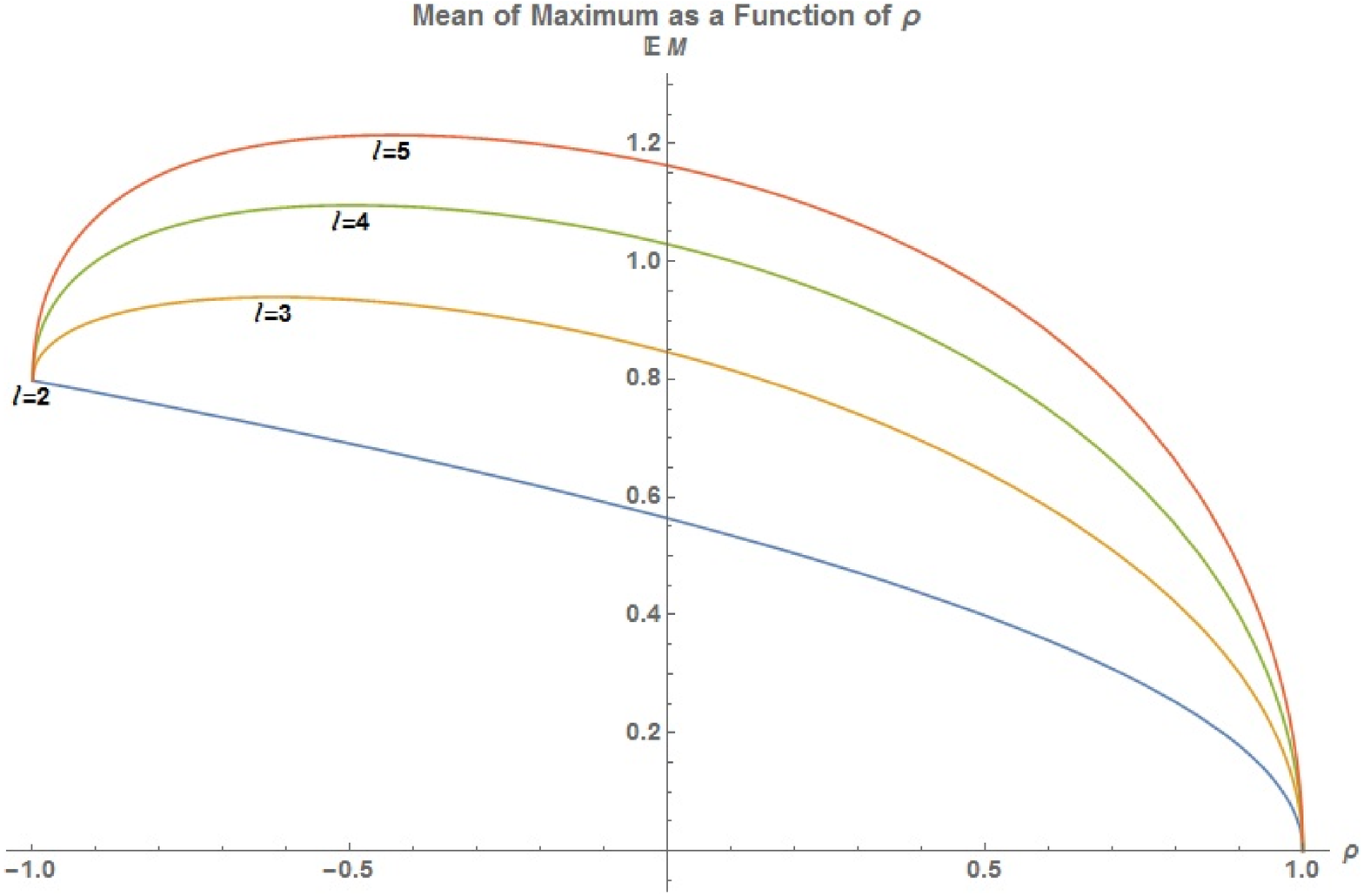}%
\caption{ Interior maximum points exist for $\ell\geq3$, but not for $\ell
=2.$}%
\end{figure}
%EndExpansion

%

%TCIMACRO{\FRAME{ftbpFU}{6.3235in}{4.0387in}{0pt}{\Qcb{Linear (with slope
%$1/\pi$ and vertical intercept $1-1/\pi$) for $\ell=2$; strictly increasing
%for $\ell\geq3$.}}{}{afonjav0.eps}{\special{ language "Scientific Word";
%type "GRAPHIC";  maintain-aspect-ratio TRUE;  display "USEDEF";
%valid_file "F";  width 6.3235in;  height 4.0387in;  depth 0pt;
%original-width 10.8923in;  original-height 6.9349in;  cropleft "0";
%filename 'afonjav0.eps';file-properties "XNPEU";}}}%
%BeginExpansion
\begin{figure}[ptb]%
\centering
\includegraphics[
height=4.0387in,
width=6.3235in
]%
{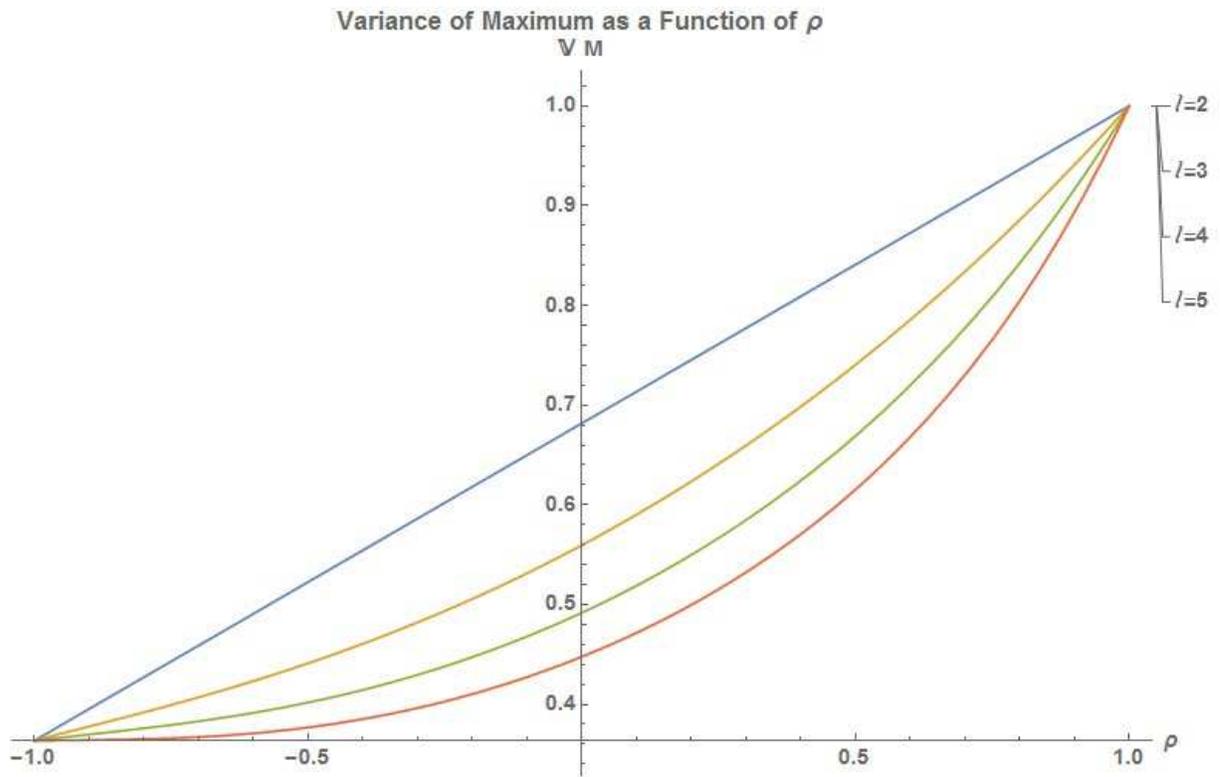}%
\caption{Linear (with slope $1/\pi$ and vertical intercept $1-1/\pi$) for
$\ell=2$; strictly increasing for $\ell\geq3$.}%
\end{figure}
%EndExpansion

\end{document}